\documentclass[12pt]{amsart}

\usepackage{amsmath,amssymb,latexsym} 
\usepackage{fullpage}
\usepackage{etoolbox}
\usepackage{enumitem}
\usepackage{color}
\usepackage[colorinlistoftodos,bordercolor=orange,backgroundcolor=orange!20,linecolor=orange,textsize=scriptsize]{todonotes}
\usepackage{enumitem}

\newtheorem{theorem}{Theorem} 
\newtheorem{corollary}{Corollary}
\newtheorem{proposition}{Proposition}

\newtheorem{lemma}{Lemma}

\theoremstyle{definition}

\theoremstyle{remark}

\newtheorem*{proof-claim}{Proof}

\newenvironment{changemargin}[2]{\begin{list}{}{%
\setlength{\topsep}{0pt}%
\setlength{\leftmargin}{0pt}%
\setlength{\rightmargin}{0pt}%
\setlength{\listparindent}{\parindent}%
\setlength{\itemindent}{\parindent}%
\setlength{\parsep}{0pt plus 1pt}%
\addtolength{\leftmargin}{#1}%
\addtolength{\rightmargin}{#2}%
}\item }{\end{list}}

\AtBeginEnvironment{proof-claim}{\vspace{-0.5cm}\footnotesize\begin{changemargin}{0.5cm}{0.5cm}}
\AtEndEnvironment{proof-claim}{\end{changemargin}}

\def\C{\mathbb{C}}

\def\C{\mathcal{C}}

 \title{An overview of complex ellipsoids }
 \author{Jorge Luis Arocha}
 \author{Javier Bracho}
 \address{Instituto de matemáticas, UNAM, Mexico}
 \author{Luis Montejano}
\address{Instituto de Matematicas, UNAM, Sede Juriquilla, Queretaro, Mexico}

\begin{document}

\maketitle 
\tableofcontents

\section{Introduction}

An ellipsoid is the image of a ball under an affine transformation. If this affine transformation is over the complex numbers, we refer to it as a \emph{complex ellipsoid}. Characterizations of  real ellipsoids have received much attention over the years, see for example de surveys \cite{Pe}, \cite{Sol}, \cite{GH} or Section 1.12 of \cite{MMO}; however, characterizations of complex ellipsoids have been studied very little. The recent interest in them begins with Gromov's proof of the  Isometric Banach Conjecture over the 
complex numbers in even dimensions \cite{G}.  Gromov proved that a complex symmetric convex body in complex space $\mathbb{C}^{2n+1}$ all whose  complex hyperplane sections through the origin are complex linearly equivalent  is a complex  ellipsoid. For that purpose, he began by proving that
the group of complex linear symmetries of a compact complex symmetric set of $\mathbb{C}^{n}$ is conjugate to a compact subgroup of $SU_n$ using the uniqueness of the complex ellipsoid of minimal volume containing a complex symmetric set.  

This paper is a review of what is known about complex ellipsoids from the point of view of convex geometry.  In particular, the proof of the Complex Banach Conjecture for dimensions $n\equiv 0,1,2$, $n\geq 5$, in \cite{BM} is outlined in Section~\ref{Sect_Banach}. This conjecture is basically a characterization of complex ellipsoids. Its proof  revealed that some ideas from the real case could be carried on to the complex case, but also that many new intriguing questions naturally arise. One of them was that if the analog of segments in the real line is discs in the complex line, the analog of convex bodies in $\mathbb{R}^n$ should be what we called \emph{bombons} in $\mathbb{C}^n$: bodies whose non-void sections with complex lines are disks. They don't give any new examples, because they turned out to be complex ellipsoids, Section~\ref{Sec:bombons} is devoted to that proof which first appeared in \cite{ABBCE}. One of the main tools used in many characterizations of complex ellipsoids is the uniqueness of minimal circumscribed complex ellipsoids of convex bodies; which is reviewed in Section~\ref{Secc:extremales}.

The study of complex ellipsoids is naturally related to the study of complex symmetry. So, characterizing and understanding complex symmetry is vital to characterizing complex ellipsoids; the next section is devoted to it. 

To ease the reading, we omit the term \emph{complex} from such concepts as lines, hyperplanes, subspaces, $k$-planes, affine and linear maps, bodies of revolution and ellipsoids when they are clearly meant to be in the context of a complex vector space. On the other hand, we will always use the term \emph{real} when we work with their counterparts over the real numbers. However,  we  still use the term complex in the titles of sections and for emphasis in definitions.

\section{Complex symmetry}

 Let $\mathbb{S}^1$ be the multiplicative group of unit complex numbers $\mathbb{C}$. 
 Let $A\subset \mathbb{C}^n$ be a set.  We say that $A$ is \emph{complex symmetric}, or simply \emph{symmetric}, 
 if there is a translated copy $A'$ of $A$ such that 
$\xi A'=A'$, for every $\xi\in \mathbb{S}^1$. In this case, if  $A'=A-x_0$, we say that $x_0$ is the center of symmetry of $A$.
If $-A$ is a translated copy of $A$, we just say that $A$ is \emph{real-symmetric}.   It will be useful to consider the empty set as a symmetric set.  Note that a \emph{ a convex body $A\subset \mathbb{C}^n$ is symmetric with center at $x_0$  if and only if for every line $L$ through $x_0$, the section $L\cap A$ is a disk centered at $x_0$.}  Of course,  the unit ball of a finite dimensional Banach space over the complex numbers $\mathbb{C}$ is symmetric. 

The following characterization of symmetry will be very useful.

\begin{theorem}\label{teosc}
A convex body $K\subset \mathbb{C}^n$  is symmetric if and only if for every $\xi\in \mathbb{S}^1$, 
$\xi K$ is a translated copy of $K$.
\end{theorem}
\proof Suppose $K$ is   symmetric, then there exists $x_0\in\mathbb{C}^n$ such that for every $\xi\in \mathbb{S}^1$, $\xi(K-x_0)=K-x_0$. This implies that 
$\xi K-\xi x_0=K-x_0$ and hence that 
$$\xi K= K + (\xi -1)x_0\,.$$ 
Thus, $\xi K$ is a translated copy of $K$. 
  
Suppose now that for every $\xi\in \mathbb{S}^1$, $\xi K$ is a translated copy of $K$. Hence 
$$\xi K= K+ \kappa_\xi\,,$$ 
where $ \kappa_\xi\in \mathbb{C}^n$  depends on $\xi$. 

For $\xi = -1$, we get that $- K= K+ \kappa_{-1}$. This implies that $K$ is real-symmetric with center $-\frac{\kappa_{-1}}{2}$. Define $x_0=-\frac{\kappa_{-1}}{2}$; it is the center of real-symmetry of $K$.

We then have that for any given $\xi\in \mathbb{S}^1$, $\xi K$  is real-symmetric with center $\xi x_0$.
On the other hand, since $\xi K= K+ \kappa_\xi$, the center of real-symmetry of $\xi K$ is $x_0 + \kappa_\xi$. Since there is a unique center of real-symmetry of $\xi K$, we have that
$$\xi x_0=x_0 + \kappa_\xi\,.$$  
Therefore, $\kappa_\xi=(\xi -1)x_0.$

To conclude the proof, note that $\xi K=K+(\xi-1)x_0,$ implies that $\xi( K-x_0)=K-x_0$ for every $\xi\in \mathbb{S}^1$ and therefore,  
$K$ is  symmetric by definition. \qed

\medskip

The analog of the following theorem but over the real numbers was originally proved by Rogers  \cite{ro1} in 1964. Using it, it is possible to prove that if all real hyperplane sections of a convex body through a point $O$ are real-symmetric, it is because the body is itself real-symmetric.  Indeed, if the point  ${O}$ is not the center of symmetry of the body, then Larman \cite{la} proved, in the real case, that the body must be a real-ellipsoid. This result  is known as the \emph{false center theorem}.  We still don't know if the false center  theorem is true over the complex numbers.  

 
\begin{theorem}[Complex Rogers' Theorem] \label{lemrg}
Let $K_1, K_2\subset \mathbb{C}^n$, $n\geq 3$, be convex bodies with ${p_1\in int(K_1)}$ and $p_2\in int(K_2)$. Suppose that for every   hyperplane $H$ through the origin,  
the section $(H+p_1)\cap K_1$ is positively homothetic to $(H+p_2)\cap K_2$. Then, $K_1$ is positively homothetic of $K_2.$ 
\end{theorem}

\proof  Asume first that $K_1$ and  $K_2$ are strictly convex bodies. 
Note that if $K_2'$ is positively homothetic to $K_2$, then $K_2'$ also satisfies the conditions of the lemma. So that using a positive homothecy, we may assume $p_1=p_2$ and  ${K_2\subset int(K_1)}$.  Now, dilate $K_2$ from $p_1$ until it touches the boundary of $K_1$. Therefore, without loss of generality we may assume that:
\begin{enumerate}
\item $K_2\subset K_1,$
\item $p_1=p_2$, and 
\item $K_1$ and $K_2$ share a common support real-hyperplane $\Delta$ at $x_0$. 
\end{enumerate}

Let $L$ be the 
 line through $x_0$ and $p_1=p_2$. Consider a hyperplane $H$ that contains $L$. By hypothesis, there is a positive homothecy  $h$ that sends $H\cap K_1$ onto $H\cap K_2$. We claim that $h$ has center $x_0$. Since such a homothecy sends support real-hyperplanes (of $H\cap K_1$) to parallel support real-hyperplanes (of $H\cap K_2$), it has to send $H\cap \Delta$ to itself because both $K_1$ and $K_2$ are strictly convex and the homothecy is positive. So, $h$ has center $x_0=K_1\cap\Delta=K_2\cap\Delta$; let $r_H>0$ be its ratio.
Since $x_0\in L$, $h$ restricts to a positive homothecy that sends $L\cap K_1$ to $L\cap K_2$. Because $L\cap int K_1\neq\emptyset$, the ratio of homothecy $r_H$ must be the same $r$ for all $H$ that contain $L$.   
Because $n\geq 3$, every point lies in a  hyperplane through $L$, so that the  homothecy with center at $x_0$ and ratio $r$ sends  $K_1$ to $K_2$ as we wished to prove. 

The proof of the non-strictly convex case requires one more step of a slightly more technical nature that follows straightforwardly from the real case given in \cite{ro1}.
\qed

\medskip 

Note that in Theorem \ref{lemrg} all homotecies are real homotecies. So, it is interesting to ask if a similar result exists for complex homotecies. Namely: \emph{ if $K_1$ and $K_2$ are convex bodies in $\mathbb{C}^n$, $n\geq 3$, containing the origin in their interior and for every hyperspace $H$, there exists a complex number $\xi$ (depending on $H$) such that
$$\xi(H\cap K_1)=H\cap K_2\,,$$
is it true that $\,\xi_0 K_1= K_2\,$, for some $\xi_0\in \C$?}

\medskip

As a corollary to Rogers' Theorem we have: 


\begin{corollary}\label{corrg}
Let $K_1, K_2\subset \mathbb{C}^n$, $n\geq 3$, be convex bodies with $p_1\in int (K_1)$ and $p_2\in int(K_2)$. Suppose that for every 
hyperspace $H$  
the section $(H+p_1)\cap K_1$ is a translated copy of ${(H+p_2)\cap K_2}$. Then, $K_1$ is a translated copy of $K_2.$ 
\end{corollary}

\proof Translations are homothecies of ratio 1. So, by Theorem \ref{lemrg}, we conclude that there is a positive homothecy $h$ sending $K_1$ to $K_2$, with
ratio of homothecy $r>0$. To prove that  $h$ is a translation we must show that $r=1$.  Let $p_3=h(p_1)$. Then, for every  hyperplane $H$ through the origin  
the section $(H+p_3)\cap K_2$ is  positively homothetic to $(H+p_2)\cap K_2$, with ratio of homothecy $r$. Let $H'$ be a  hyperplane through the origin with the property that $(H'+p_3)=(H'+p_2)$. Hence, $(H'+p_3)\cap K_2=(H'+p_2)\cap K_2$, and being positively homothetic with ratio of homothecy $r>0$, we conclude that $r=1$. \qed

\begin{theorem}\label{lemBB}
In dimension $n \geq 3$, a convex body $K\subset \mathbb{C}^n$ all of whose  hyperplane sections through an interior point $x_0$ are  symmetric, 
is itself  symmetric. 
\end{theorem}

\proof  To prove that $K$ is   symmetric, by Theorem \ref{teosc} it is enough  to prove that for every $\xi\in \mathbb{S}^1$, $\xi K$ is a translated copy of $K$. 

Fix  $\xi\in \mathbb{S}^1$, thus $\xi x_0$ is an interior point of $\xi K$.  It is enough to prove that for every  hyperspace $H$, $(H+x_0)\cap K$ is a translated copy of  $(H+\xi x_0) \cap \xi K$; because this implies, by Corollary \ref{corrg}, that $\xi K$ is a translated copy of $K$.

By hypothesis, 
$(H+x_0)\cap K$ is  symmetric, hence by Theorem \ref{teosc}, $(H+x_0)\cap K$ is a translated copy of $\xi\big((H+x_0)\cap K\big)$, but 
$$\xi\big((H+x_0)\cap K\big)= (\xi H+\xi x_0)\cap \xi K= (H+\xi x_0)\cap \xi K\,.$$ \qed

We now turn our attention to characterizing complex symmetry by means of projections.

 
\subsection{Complex symmetry and projections}

The  affine image of a  symmetric set is a  symmetric set. Essentially this is so because, if $\xi K=K$, for every $\xi\in \mathbb S^1$, and $f$ is a linear map then $\xi f(K)=f(\xi K)=f(K)$. In particular, the orthogonal projection of a symmetric set onto a  line is always a disk. The following simple lemma will be useful in the sequel.

\begin{lemma}\label{lemSD}
A convex body $K\subset \mathbb{C}^n$ is  symmetric if and only if it has a translated copy all whose orthogonal projections onto $1$-dimensional  subspaces are disks centered at the origin.
\end{lemma}

\proof The necessity is clear. For sufficiency, assume $\pi(K)$ is a disk centered at the origin for every orthogonal projection $\pi$ onto a $1$-dimensional subspace.

We shall prove that for every $\xi\in \mathbb S^1$, $\xi K=K$.

Suppose not. Then there is $x\in K$ and $\xi\in \mathbb S^1$ such that 
$\xi x\notin K$. By convexity there is a 
 real-hyperplane $\Delta$ in $\mathbb R^{2n}=\mathbb C^n$ that separates $\xi x$ from   
$K$. Let $H$ by the unique  hyperplane contained in $\Delta$. Let $\pi:\mathbb C^n \to L$ be the orthogonal projection to the $1$-dimensional  subspace $L$ orthogonal to $H$. Then the real line $\pi(\Delta)$ separates $\pi(K)$ from $\pi(\xi x)=\xi(\pi(x))$ in the  line $L$. On the other hand, $\pi(x)$ lies in the disk
$\pi(K)$ centered at the origin, therefore $\xi(\pi(x))$ lies in $\pi(K)$, which is a contradiction. \qed

\begin{theorem}\label{teoSim} 
For $k\geq 2$, a convex body all whose orthogonal projections onto  $k$-planes are  symmetric is  symmetric.
\end{theorem}

\proof Consider a convex body $K\subset \mathbb C^n$, $n\geq 3$, all whose orthogonal projections onto hyperplanes are  symmetric. For every line $L$ through the origin, let $c_L$ be the center of the orthogonal projection of $K$ to $L^\perp$, the orthogonal hyperplane to $L$, and let $L^\prime$ be the parallel line to $L$ through $c_L$. We claim that
$$\bigcap_L L^\prime \text{ is a single point } \{x_0\}.$$

Let $L_1$ and $L_2$ be two different lines through the origin. The orthogonal projection of $K$ to $L_1^\perp\cap L_2^\perp$ is symmetric because it is the affine image of symmetric convex bodies (in $L_1^\perp$ and $L_2^\perp$); let $a$ be its center of symmetry and
let $\Gamma$ be the orthogonal plane to $L_1^\perp\cap L_2^\perp$ passing through $a$. The lines $L_1^\prime$ and $L_2^\prime$ are both in $\Gamma$ because the symmetry centers $c_{L_1}$ and $c_{L_2}$ are orthogonally projected to $a$. Therefore, since $L_1^\prime$ and $L_2^\prime$ are not parallel, they intersect in a point $\{x_0\}=L_1^\prime\cap L_2^\prime$.

Consider a third  line $L_3$ through the origin, linearly independent of $L_1$ and $L_2$; it exists because $n\geq 3$. By the above argument, $L_3^\prime$ intersects both $L_1^\prime$ and $L_2^\prime$. But by  linear independence, $L_3^\prime$ intersects $\Gamma$ in at most one point, so we must have that $x_0\in L_3^\prime$.

Finally, for the general line $L$ through the origin, $L$ is linearly independent of at least one pair of the lines $L_1, L_2, L_3$. So that the preceding argument yields that $x_0\in L^\prime$. This proves that $\bigcap_L L^\prime =\{x_0\}$ as we claimed.

Suppose without loss of generality that $x_0$ is the origin. Then, for every $(n-1)$-plane $H$ through the origin, the orthogonal projection of $K$ 
onto $H$ is symmetric with center at the origin. This immediately implies that  orthogonal projections of $K$ onto $1$-dimensional subspaces are disks centered at the origin and therefore, by Lemma \ref{lemSD}, $K$ is symmetric. 

The theorem now follows by induction, because a convex body all whose orthogonal projections onto $k$-planes are symmetric has the property that all orthogonal projections onto $(k+1)$-planes are symmetric, and so on. 
\qed

At this point the following two questions arise.  
\emph{ Is a convex body all whose orthogonal projections onto lines are disks,  symmetric?}  
The projection of an ellipsoid from an $(n-2)$-affine plane on a line is a disk. Is the converse true? That is, \emph{ if all projections of a convex body on  lines are disks, is the body an ellipsoid?}

\subsection{Complex symmetry and complex ellipsoids}
\begin{lemma}\label{symellii}
A symmetric real ellipsoid is an ellipsoid.
\end{lemma}

\proof We need to recall some  facts about real ellipsoids, within the real context.  
Let $E\subset \mathbb{R}^n$ be a $n$-dimensional real ellipsoid centred at the origin. For every $k$-dimensional real subspace, $H\subset \mathbb{R}^n$ with $1\leq k< n$, there exists a complementary $(n-k)$-dimensional real subspace, $L$ of $\mathbb{R}^n$, called its {\em polar subspace with respect to $E$}, such that  
$$\partial E \cap L=\{x\in\mathbb{R}^n \mid H+x \mbox{ is a k-dimensional real plane tangent to } \partial E \mbox{ at } x\}$$ 
(this set is called the {\em shadow boundary of  $E$ in the direction $H$}). Moreover, $H$ is the polar subspace of $L$ with respect to $E$, and the section $L\cap E$ is a $(n-k)$-dimensional real ellipsoid with the following property:  for every $(n-k)$-plane $L'$, parallel to $L$, the corresponding section $L'\cap E$ is either the empty set, a point in $H$ or a real ellipsoid homothetic to $L\cap E$ and centred at $H$.  For more about shadow boundaries see Section 1.12.2  of \cite{MMO}. 
 
 Clearly, every ellipsoid is a real ellipsoid which is symmetric. 
Let ${K\subset \mathbb{C}^n=\mathbb{R}^{2n}}$ be a  symmetric real ellipsoid centered at the origin. By induction on the complex dimension $n$, 
we will prove that there is a linear isomorphism $g\in GL(n,\mathbb{C})$
such that $g(K)$ is a ball. The statement is true for $n=1$. Suppose it is true for dimension $n-1$, we shall prove it for dimension $n$.  

Assume the diameter of $K$ is $h$, and let $[-u,u]$ be a diameter of $K$; let $L$ be the unique  line containing the vector $u$. By hypothesis $D=L\cap K$ is a disk centred at the origin all of whose diameters are also diameters of $K$. This implies that the polar to $L$ with respect to $E$ is the  hyperplane, $H$, orthogonal to $L$. Then, for every affine  line $L'$ orthogonal to $H$ and touching $int(K)$, the section $L'\cap K$ is a disk with center at $H$.  

By induction we have that $H\cap K$ is an ellipsoid. Therefore, using a linear isomorphism, we may assume that $H\cap K$ is a $(2n-2)$-dimensional ball of diameter $h$.  To conclude the proof of the lemma, we prove that $K$ is a ball. 

Let $\ell$ be a real line subspace contained in $H$ and let $\Delta$ be the $3$-dimensional real subspace generated by $\ell$ and $L$.  Since $(L+x)\cap K$ is a disk with centre at $\ell$ for every $x\in \ell\cap int(K) $, 
$\Delta \cap K$ is a real ellipsoid of revolution with axis the line $\ell$.  Since the three axes of this ellipsoid are equal, this implies,  that $\Delta\cap K$ is a $3$-dimensional ball with centre at the origin. Since this holds for every real $3$-dimensional real subspace containing $L$, we have that $K$ is a ball, as we wished. 
 \qed

 In 1955, in his book; The Geometry of Geodesics \cite{Bus}, Busemann  proved the following Kubota's Theorem: \emph{if all real planes passing through a point $O$ intersect a convex body (in $\mathbb{R}^3$) in real ellipses then the body is a real ellipsoid}. We will start proving the analog of this theorem for ellipsoids with the additional hypothesis that the point
$O$ is the  center of symmetry of the body. Later, in Theorem~\ref{teoksec}, we will get rid of this auxiliary hypothesis.

 \begin{proposition}\label{prop:Celipsoo} 
If $K\subset \mathbb{C}^{n+1}$, $n\geq 2$, is a  symmetric convex  body all whose  hyperplane sections through the center are  ellipsoids, then $K$  is an  ellipsoid. 
\end{proposition}

\proof By Lemma \ref{symellii}, it is enough to prove that $K$ is a real ellipsoid.  
Consider $\mathbb{C}^n = \mathbb{R}^{2n}$. By Theorem 2.12.4 of \cite{MMO}, it is enough to prove that every real two dimensional subspace intersects $K$ in an ellipse.  Let $\Pi$ be a two dimensional real plane generated by $\{v_1, v_2\}$.  If $\Pi$ is a  line, $\Pi\cap K$ is a ball, so assume it is not. Let $L_i$ be the  line containing $v_i$, $i=1,2$. Consequently,  $\Pi$ is contained in the  plane $P$ generated by $\{L_1, L_2\}$. By hypothesis,  $P\cap K$ is a section of an ellipsoid and hence is itself an ellipsoid.
This implies that $\Pi\cap K$ is an ellipse. Therefore, $K$ is a real ellipsoid. 
 \qed

\section{Extremal inscribed and circumscribed complex ellipsoids}\label{Secc:extremales}
 Let $A$ be a non-flat, compact subset of  $\mathbb{C}^{n}$.  An ellipsoid $%
\mathfrak{E}$ is called \emph{circumscribed} if $A\subset \mathfrak{E}$, and
it is  \emph{inscribed} if $\mathfrak{E}\subset \widehat{A}$, where $\widehat{A}$ denotes the convex closure of $A$. We
say that $\mathfrak{E}$ is a \emph{minimal circumscribed ellipsoid} (MiCE) if it has the minimal volume 
among all circumscribed ellipsoids. On the
other hand, $\mathfrak{E}$ is a \emph{maximal inscribed ellipsoid} (MaIE) if $\mathfrak{E}$ has maximal volume among all inscribed ellipsoids.

The purpose of this section is to prove that if a convex set in $\mathbb{C}^{n}$ contains two
inscribed ellipsoids of maximal volume, then one is a translate of the
other. On the other hand, the circumscribed ellipsoid of minimal
volume is unique. The first proofs of uniqueness of MiCE and MaIE for real ellipsoids in its full generality seem to have appeared
independently in \cite{DLL} and \cite{Zaguskin} and for complex ellipsoids in \cite{ABBE}.  For real ellipsoids MiCE and MaIE are known
as L\"{o}wner--John Ellipsoids and have many applications in several areas
of mathematics (see, e.g., \cite{HM}, \cite{Pe}, Section~2.12 of \cite{MMO} and \cite{H} with the references therein).

The existence of these ellipsoids follows from standard arguments. There is a ball big enough to contain $A$ and a non-zero sphere contained in $\widehat{A}$.
 All ellipsoids
are easily parametrized by a matrix and a vector; among them, we can
consider only those that are contained in the
big ball. This set is compact in the parameter space. Moreover, the volume
function is continuous so that the existence of a MiCE and a MaIE follows.

For $x=\left( x_{1},...,x_{n}\right) \in \mathbb{C}^{n}$, let $%
\overline{x}=\left( \bar{x}_{1},...,\bar{x}_{n}\right) $, where the
bar above denotes the complex conjugate. Also, we denote by $x\odot y$ the Hadamar product of $x$ and $y$; that is, the coordinatewise product (see for example \cite{Horn} Chapter 5). 

A \emph{scalar product} in $\mathbb{C}^{n}$ is a sesquilinear, Hermitian, positive-definite functional $\mathbb{C}^{n}\times \mathbb{C}^{n}\rightarrow \mathbb{C}$ denoted by $\left\langle x\cdot y\right\rangle $ for
any $x$ and $y$ in $\mathbb{C}^{n}$. For each scalar
product there is a matrix $M$ (Hermitian, positive-definite) such that $%
\left\langle x\cdot y\right\rangle = x ^{T}M%
\overline{y}$. In the case that $M$ is the identity matrix, the
scalar product is the usual Hermite's product in $\mathbb{C}^{n}$. It is well known that the eigenvalues of $M$ are real positive numbers and
so is its determinant. It is also known that $M$ is diagonalizable by a
unitary transformation.

An \emph{ellipsoid} (centered at the origin) is a set%
\begin{equation}
\left\{  x \in \mathbb{C}^{n}\mid   x^{T}M\overline{ x }%
\leq 1\right\} \,,  \label{Elipsoid}
\end{equation}%
where $M$ is Hermitian and positive-definite.
If $M$ is the identity matrix, this ellipsoid is the \emph{unit
ball} $\mathfrak{B}$.

The set of unitary
transformations is the subgroup of $GL\left( \mathbb{C}^{n}\right) $ that
preserves the unit sphere. The modulus of a scalar $\xi \in \mathbb{C}$
will be denoted by $\left\vert \xi \right\vert $.  We will denote by $\left\Vert \cdot \right\Vert 
$ the usual norm in $\mathbb{C}^{n}$, that is, the norm defined by Hermite's
product.

Let $B\overset{\text{def}}{=}\sqrt{M^{T}}$. We have $M^{T}=BB$ and
hence, $M=B^{T}B^{T}=B^{T}\overline{B}$. From this, we obtain
\begin{equation*}
x^{T}M\overline{x}=x^{T}B^{T}\overline{B}%
\overline{x}=\left( Bx\right) ^{T}\overline{\left( B%
x\right) }=\left\Vert Bx\right\Vert ^{2}\,.
\end{equation*}%
Therefore, the ellipsoids in $\mathbb{C}^{n}$ can also be written in the form%
\begin{equation*}
\{ B^{-1}u\mid u\in \mathfrak{B} \}.
\end{equation*}

If we use a unitary transformation to bring $M$ to the diagonal form, then
this can be rewritten as%
\begin{equation}
\mathfrak{El}\left( \lambda \right) \overset{\text{def}}{=}\left\{ 
\lambda \odot u\mid u\in \mathfrak{B} \right\}\,,   \label{Forma parametrica}
\end{equation}
where $\mathbf{\lambda }$ is a vector in $\mathbb{R}_{+}^{n}$. Since $M^{{\frac12}}=B$, the two forms \ref{Elipsoid} and \ref{Forma parametrica} are related
by the fact that $\lambda $ is the diagonal of $M^{-{\frac12}}$. The map 
$x\mapsto \lambda \odot x=M^{-{\frac12}}x$ is an invertible linear map in $GL\left( \mathbb{C}^{n}\right)$
which maps the unit ball $\mathfrak{B}$ into
the ellipsoid $\mathfrak{El}\left(\lambda \right)$. Therefore 
\begin{equation}
Vol\mathfrak{El}\left( \lambda \right) =\det M^{-{\frac12}}
Vol\mathfrak{B} \,.  \label{Volumen-det}
\end{equation}

The results that follow do not depend on the translation-invariant measure chosen to define the volume. We just need the validity of equation~\ref{Volumen-det}.

Let $\det\lambda $ denote the product of coordinates of $\lambda $. Of course, we have ${\det \lambda =\det A^{-{\frac12}}}$ and%
\begin{equation*}
Vol\mathfrak{El}\left( \lambda \right) =\det \lambda\, 
Vol\mathfrak{B} .
\end{equation*}

Translates of ellipsoids centered at the origin are also  \emph{ellipsoids}, and translations do not change volume.

\subsection{Auxiliary lemmas}
For the proofs of Theorems \ref{MaIE} and \ref{MiCE}, we need first  two lemmas whose proofs are
straightforward computations. Let $\mathbf{1}$ be the vector which
has all coordinates equal to $1$.

\begin{lemma}
\label{volumen}If $\mathbf{\lambda }\in \mathbb{R}_{+}^{n}$ is such that
$\det \mathbf{\lambda }=1$ and $\mathbf{\lambda }\neq \mathbf{1}$, then%
\begin{equation*}
\det \left( \frac{\mathbf{\lambda }+\mathbf{1}}{2}\right) >1\text{.}
\end{equation*}
\end{lemma}

\proof
We have $\left( \lambda _{i}-1\right) ^{2}\geq 0$ and therefore $\left(
\lambda _{i}+1\right) ^{2}\geq 4\lambda _{i}$. With equality only when $%
\lambda _{i}=1$. Taking the product, we get 
$$\det \left( \mathbf{\lambda }+%
\mathbf{1}\right) ^{2}>4^{n}\det \mathbf{\lambda }$$ 
which is the same as 
$$
\det \left( \mathbf{\lambda }+\mathbf{1}\right) >2^{n}\sqrt{\det \mathbf{%
\lambda }}=2^{n}\,.$$ 
The lemma follows.
\qed

\begin{lemma}
\label{Lema para el afin}Let $c,x$ be complex numbers and $\lambda $ be a real
number. Then,%
\begin{equation*}
\lambda \left\vert x\right\vert ^{2}+\left\vert x-c\right\vert ^{2}=\left(
\lambda +1\right) \left\vert x-\frac{c}{\left( \lambda +1\right) }%
\right\vert ^{2}+\left( \frac{\lambda }{\lambda +1}\right) \left\vert
c\right\vert ^{2}\,.
\end{equation*}
\end{lemma}

\proof
We have%
\begin{equation*}
\begin{array}{c}
\lambda \left\vert x\right\vert ^{2}+\left\vert x-c\right\vert ^{2}=\lambda
\left\vert x\right\vert ^{2}+\left( x-c\right) \overline{\left( x-c\right) }
\\ 
=\left( \lambda +1\right) \left\vert x\right\vert ^{2}-\left( x\overline{c}+c%
\overline{x}\right) +\left\vert c\right\vert ^{2}%
\end{array}%
\end{equation*}%
and completing the square, we obtain%
\begin{equation*}
\begin{array}{l}
=\left( \lambda +1\right) \left( \left\vert x\right\vert ^{2}-\dfrac{\left( x%
\overline{c}+c\overline{x}\right) }{\left( \lambda +1\right) }+\left\vert 
\dfrac{c}{\left( \lambda +1\right) }\right\vert ^{2}-\left\vert \dfrac{c}{%
\left( \lambda +1\right) }\right\vert ^{2}\right) +\left\vert c\right\vert
^{2} \\ 
=\left( \lambda +1\right) \left\vert x-\dfrac{c}{\left( \lambda +1\right) }%
\right\vert ^{2}+\left( \dfrac{\lambda }{\lambda +1}\right) \left\vert
c\right\vert ^{2}.%
\end{array}%
\end{equation*}
\qed

\subsection{Maximal inscribed complex ellipsoids.}

\begin{theorem}\label{MaIE}
Let $A$ be a non-flat compact in $\mathbb{C}^{n}$. Let $\mathfrak{E}_{1}$ and $\mathfrak{E}_{2}$ be two MaIE contained in $\widehat{A}$. Then,
there is a vector $c\in \mathbb{C}^{n}$ such that $\mathfrak{E}_{2}=%
\mathfrak{E}_{1}+c$.
\end{theorem}

\proof
Using a suitable affine transformation we can assume that $\mathfrak{E}_{1}$
is the unit ball $\mathfrak{El}\left( \mathbf{1}\right)$. Let $c\in \mathbb{C}^{n}$ be the center of $\mathfrak{E}_{2}$.
We can use a unitary transformation to diagonalize the matrix of $\mathfrak{E}_{2}-c$. And therefore, 
$\mathfrak{E}_{2}=\mathfrak{El}\left( 
\lambda \right) +c$ for some $\lambda \in\mathbb{R}_{+}^{n}$.

Let us first prove that the ellipsoid $\mathfrak{E}_{3}=\mathfrak{El}\left({\frac12}\left(\lambda +\mathbf{1}\right) \right) +{\frac12}c$ is contained in the convex closure of $\mathfrak{E}_{1}\cup 
\mathfrak{E}_{2}$, that is, 
\begin{equation*}
\mathfrak{E}_{3}=\mathfrak{El}\left( \frac{\lambda +\mathbf{1}}{2}\right) +\frac{%
c}{2}\subset \widehat{\mathfrak{El}\left( \mathbf{1}\right) \cup
\left( \mathfrak{El}\left( \lambda \right) +c\right) }\,.
\end{equation*}%

For any $x\in \mathfrak{E}_{3}$ there is $u\in 
\mathfrak{B} =\mathfrak{El}\left( \mathbf{1}\right) =\mathfrak{E}_{1}$, such that
\begin{equation*}
x=\frac{\lambda +\mathbf{1}}{2}\odot u+\frac{c}{2}\,.
\end{equation*}

The point $y=\lambda \odot u+c$ is in $\mathfrak{E}_{2}$ and
\begin{equation*}
\frac{y+u}{2}=\frac{\lambda \odot u+c+u}{2}=x\,.
\end{equation*}
This means that $x$ is the middle point of the segment joining
 $y$ and $u$, and proves that $\mathfrak{E}_{3}\subset \widehat{\mathfrak{E}_{1}\cup \mathfrak{E}_{2}}$, as we wished. 

Since $\mathfrak{E}_{1}\cup \mathfrak{E}_{2}\subset \widehat{A}$, we have 
$\mathfrak{E}_{3}\subset \widehat{\mathfrak{E}_{1}\cup \mathfrak{E}_{2}}\subset 
\widehat{A}$. 

We know that $Vol\left( \mathfrak{E}_{1}\right) =Vol\left( 
\mathfrak{E}_{2}\right) =\det \lambda Vol\left( \mathfrak{E}_{1}\right)$. 
Therefore, $\det \lambda =1$. 
On the other hand, 
$Vol\left( \mathfrak{E}_{3}\right) =
\det \left({\frac12}\left( \lambda +
\mathbf{1}\right) \right) Vol\left(\mathfrak{E}_{1}\right)$. 
If $\lambda \neq \mathbf{1}$, 
Lemma~\ref{volumen}  implies that 
$\det \left({\frac12}\left( \lambda +\mathbf{1}\right) \right) >1$, 
and hence 
$Vol\left( \mathfrak{E}_{3}\right) > Vol\left( \mathfrak{E}_{1}\right)$; 
but this contradicts that $\mathfrak{E}_{1}$ is a MaIE. So, we have 
$\lambda =\mathbf{1}$ and therefore, $\mathfrak{E}_{2}=\mathfrak{E}
_{1}+c$, which proves our theorem.
\qed

Now, we shall see that symmetry guarantees uniqueness of the MaIE. 

\begin{theorem}
If $A$ is a non-flat, compact symmetric set in $\mathbb{C}^{n}$, then its
MaIE is unique.
\end{theorem}

\proof
We can assume that the center of $A$ is the origin. Let $\mathfrak{E}$ be
the unique real MaIE in $A$. Let $\zeta $ be a scalar
of modulus $1$. Since ${\zeta \mathfrak{E}\subset \zeta A=A}$ and 
$Vol\left( \zeta \mathfrak{E}\right) =Vol\left( \mathfrak{E}\right)$
then, the uniqueness of $\mathfrak{E}$ implies that $\mathfrak{E}$ is
symmetric. Using Lemma \ref{symellii}, we get that $\mathfrak{E}$ is not only a real ellipsoid but a complex one.
There is not another  MaIE because any ellipsoid is also a real ellipsoid.
\qed

\subsection{Minimal circumscribed complex ellipsoids}

The uniqueness of a minimal  ellipsoid containing a symmetric set was first proved by Gromov \cite{G}, when he was proving that the group of  symmetries of a compact symmetric set of $\mathbb{C}^{n}$ is conjugate to a compact subgroup of $SU_n$.
Here we prove the uniqueness without the additional hypothesis of symmetry.

\begin{theorem}\label{MiCE}
Let $A$ be a non-flat, compact set in $\mathbb{C}^{n}$. Let $%
\mathfrak{E}_{1}$ and $\mathfrak{E}_{2}$ be two MiCE containing $A$. Then, 
$\mathfrak{E}_{1}=\mathfrak{E}_{2}$.
\end{theorem}

\proof
Using an affine transformation, we can make $\mathfrak{E}_{2}$ the unit
ball. Then, using a unitary
transformation, we can diagonalize the matrix of $\mathfrak{E}_{1}$. Finally, we translate so that the center of $\mathfrak{E}_{1}$ is the origin and the center of $\mathfrak{E}_{2}$ is some vector $c\in \mathbb{C}^{n}$.
Then,
\begin{equation*}
\begin{array}{l}
\mathfrak{E}_{1}=\left\{ x\in \mathbb{C}^{n}\mid \sum \lambda _{i}\left\vert
x_{i}\right\vert ^{2}\leq 1\right\} =\mathfrak{El}\left( \beta \right)  \\ 
\mathfrak{E}_{2}=\left\{ x\in \mathbb{C}^{n}\mid \sum \left\vert
x_{i}-c_{i}\right\vert ^{2}\leq 1\right\} =\mathfrak{El}\left( \mathbf{1}\right) +c\,,%
\end{array}%
\end{equation*}%
where $\mathbf{\beta }=\mathbf{\lambda }^{-{\frac12}}$. Therefore,
\begin{equation*}
A\subset \mathfrak{E}_{1}\cap \mathfrak{E}_{2}\subset \mathfrak{E}_{3}%
\overset{\text{def}}{=}\left\{ x\in \mathbb{C}^{n}\mid \sum \lambda
_{i}\left\vert x_{i}\right\vert ^{2}+\left\vert x_{i}-c_{i}\right\vert
^{2}\leq 2\right\} .
\end{equation*}%
Using Lemma \ref{Lema para el afin}, the inequality is transformed into%
\begin{equation}
\sum \left( \lambda _{i}+1\right) \left\vert \left( x_{i}-\frac{c_{i}}{%
\left( \lambda _{i}+1\right) }\right) \right\vert ^{2}\leq 2-\sum \left( 
\frac{\lambda _{i}}{\lambda _{i}+1}\right) \left\vert c_{i}\right\vert
^{2}\leq 2 \,.  \label{desigualdad}
\end{equation}

Let%
\begin{equation*}
\mathfrak{E}_{4}\overset{\text{def}}{=}\left\{ x\in \mathbb{C}^{n}\mid \sum
\left( \lambda _{i}+1\right) \left\vert \left( x_{i}-\frac{c_{i}}{\left(
\lambda _{i}+1\right) }\right) \right\vert ^{2}\leq 2\right\} \,,
\end{equation*}%
which is an ellipsoid. From inequality \ref{desigualdad}, we obtain that 
\begin{equation*}
A\subset \mathfrak{E}_{1}\cap \mathfrak{E}_{2}\subset \mathfrak{E}%
_{3}\subset \mathfrak{E}_{4} \,.
\end{equation*}%
Let $\mathbf{z}\in \mathbb{C}^{n}$ be the vector with coordinates 
$c_{i}/\left( \lambda _{i}+1\right) $. We have 
\begin{equation*}
\mathfrak{E}_{4}=\mathbf{z}+\left\{ x\in \mathbb{C}^{n}\mid \sum \left(
\lambda _{i}+1\right) \left\vert x_{i}\right\vert ^{2}\leq 2\right\} =
\mathfrak{El}\left( \left( \frac{\mathbf{\lambda }+1}{2}\right) ^{-{\frac12}}\right) +\mathbf{z}\,.
\end{equation*}
Since $\det \mathbf{\beta }=1$, we obtain $\det \mathbf{\lambda }=1$. If 
$\mathbf{\lambda }\neq \mathbf{1}$, the hypothesis of Lemma \ref{volumen} are
satisfied, then
\begin{equation*}
\Delta \overset{\text{def}}{=}\det \left( \frac{\mathbf{\lambda }+1}{2}%
\right) >1\,.
\end{equation*}%
Therefore $Vol\mathfrak{E}_{4}=\Delta ^{-{\frac12}} Vol\mathfrak{E}_{2}< Vol\mathfrak{E}_{2}$, which contradicts
the minimality of $\mathfrak{E}_{2}$.

So, we conclude that $\lambda =\mathbf{1}$, $\mathfrak{E}_{1}$ is
the unit ball and $\mathfrak{E}_{2}=\mathfrak{E}_{1}+c$. In
this case, the inequality (\ref{desigualdad}) transforms to the following%
\begin{equation*}
\sum \left\vert \left( x_{i}-\frac{c_{i}}{2}\right) \right\vert ^{2}\leq
1-\sum \left\vert \frac{c_{i}}{2}\right\vert ^{2}.
\end{equation*}%
Therefore, $\mathfrak{E}_{3}$ is a ball with center in $c/2$, and
radius 
\begin{equation*}
r\overset{\text{def}}{=}\sqrt{1-\sum \left\vert \frac{c_{i}}{2}\right\vert
^{2}}\,.
\end{equation*}%
If $c$ is not the origin then $r<1$ and $\mathfrak{E}_{3}$ has volume strictly less
than that of the unit ball $\mathfrak{E}_{1}$. Therefore, 
$c=\mathbf{0}$ and $\mathfrak{E}_{2}=\mathfrak{E}_{1}.$
\qed

\section{ A  characterization of complex ellipsoids}\label{Sec:bombons} 

The following characterization of the  ellipsoids has no analogue over the real numbers. It was unexpected for us, and was first proved in \cite {ABBCE}.

 \smallskip 
 
 \begin{theorem}\label{teoBom}
A convex body $K\subset \mathbb{C}^n$ is an ellipsoid  if and only if for every  line $L\subset \mathbb{C}^n$, $L\cap K$ is either empty, a point or a disk. 
\end{theorem}

For expository reasons, it is convenient to define that a \emph{ bombon} is a convex body $K\subset \mathbb{C}^n$ such that for every line $L\subset \mathbb{C}^n$, $L\cap K$ is either empty or a disk; where points are regarded as disks of radius zero. Then, the theorem simply says that bombons and ellipsoids are one and the same.

Note that $1$-dimensional ellipsoids and $1$-dimensional  bombons are disks. Furthermore, ellipsoids are bombons because balls are bombons and the image of a bombon under an affine map is a bombon. The basic reason for this fact is that affine maps between lines send disks to disks.

The proof of Theorem~\ref{teoBom} is by induction. The general case will follow from Proposition \ref{prop:Celipsoo}, and the first case, $n=2$, takes most of the work using ideas of its own.

An \emph{abstract linear space} consists of a set $X$  together with a distinguished family of subsets, called \emph{lines}, satisfying the following property:  given different $x,y\in X$, there is a unique line $L$ containing $x$ and  $y$.  Typical examples of abstract linear spaces are the euclidean $n$-space, the complex $n$-space and their corresponding projective spaces; with their classic lines as lines. 

Our interest in these objects lies in the abstract linear space $LS^3$, where $X=\mathbb S^3$, the unit sphere of $\mathbb C^2$, and a line is the intersection with $\mathbb S^3$ of a line of  $\mathbb C^2$, when it is has more than a single point; which is a flat circle in $\mathbb S^3$.  It is an abstract linear space because two different points in $\mathbb S^3$ lie in a unique line.

A subset $A$ of an abstract linear space $X$ is \emph{linearly closed} if for any different $x,y\in A$, the line through $x$ and $y$ is contained in $A$.  Since the intersection of linearly closed subsets is linearly closed, given any $Y\subset X$, there exists a unique minimal linearly closed subset containing $Y$, called its \emph{linear closure}.

\begin{lemma}\label{lemclo}
Suppose $K\subset\mathbb{C}^2$ is a bombon with the property that the ellipsoid of minimal volume containing it is the unit ball. Then $K\cap \mathbb S^3$ is a linearly closed subset of the abstract linear space $LS^3$.
\end{lemma}
\proof 
Let $x,y\in K\cap \mathbb S^3$, $x\neq y$, and let $L$ be the  line through $x$ and $y$. We have to prove that $L\cap \mathbb S^3\subset K$. 

Recall that $\mathfrak{B}$ is the unit ball in $\mathbb{C}^2$, whose boundary is $\mathbb S^3$. By hypothesis, $L\cap K$ and $L\cap \mathfrak{B}$ are two disks in $L$ which share two different points, $x$ and $y$ on their boundary. Then,  
$ K \subset \mathfrak{B}$ implies that $L\cap K = L\cap \mathfrak{B}$, and therefore that $L\cap \mathbb S^3\subset K$. \qed

\medskip


\begin{lemma}\label{prop:lin_cl}
The linearly closed proper subsets of the abstract linear space $LS^3$ are single points and lines.
\end{lemma}

The proof of this lemma follows immediately from the next two lemmas. 

We need to make precise some standard definitions. 
A topological space $U$ can be topologically embedded in $\mathbb C$, if there is a continuous injective map $f: U\to\mathbb C$. A topological space $A$ is locally embedded in $\mathbb C$ if for every point $x\in A$, there is a neighborhood $U$ of $x$ in $A$ which can be topologically embedded in $\mathbb C$.

\begin{lemma}
Let $A$ be a linearly closed proper subset of the abstract linear space $LS^3$. Then $A$ is embedded in $\mathbb C$.
\end{lemma}

\proof Since $A$ is proper, there exists a point $w \in\mathbb S^3$ such that $w\notin A$.  Let $L$ be the line through $w$ and the origin and let  $L^\perp$ be the orthogonal line to $L$ at the origin, which we identify with $\mathbb C$. 

Let $\pi:\mathbb S^3\setminus \{w\} \to L^\perp$  be the geometric projection from $w$, that is, given $x\in \mathbb S^3$ different  from $w$, $\pi(x)$ is the intersection with $L^\perp$ of the line through $x$ and $w$; well defined because this  line is not parallel to  $L^\perp$.

We claim that $\pi |_A$ is an embedding. Indeed, if $x, y \in A$ are such that $\pi(x)=\pi(y)$ then their lines through $w$ coincide. If $x$ and $y$  were different points, it would imply that $w\in A$ because $A$ is linearly closed, which contradicts the choice of $w$.        \qed 

\begin{lemma}
Let $A$ be a linearly closed non-empty subset of the abstract linear space $LS^3$ which is not a point or an line. Then $A$ is not locally embedded in $\mathbb C$.
\end{lemma}

\proof $A$ has at least two points, so it contains an line $\ell\subset A\subset\mathbb S^3$. But since $A$ is not an line, there is yet another point $a\in A\setminus\ell$.

For every $x\in \ell$ let $\ell_x$ be the line through $a$ and $x$.  Note that $\ell_x\cap \ell=\{x\}$, and if $x\not=x'\in \ell$,
then $\ell_x\cap\ell_{x'}=\{a\}$, because two lines intersect at most in one point. Let 
$$\Omega=\bigcup_{x\in \ell}\ell_x \subset A\subset \mathbb S^3\,.$$  
Since every abstract line of $LS^3$ is actually a circle, then $\Omega\setminus\{a\}=\bigcup_{x\in \ell}(\ell_x\setminus\{a\})$ is an $\mathbb R$-bundle over $\mathbb S^1$. This implies that, either $\Omega$ is homeomorphic to the real projective plane or $\Omega$ is homeomorphic to a closed cylinder modulo its boundary,
$\frac{\mathbb S^1\times [0,1]}{\mathbb S^1\times\{0,1\}}$, where the point $a\in \Omega$ is precisely $\frac{\mathbb S^1\times \{0,1\}}{\mathbb S^1\times\{0,1\}}$.

Since $\Omega\subset \mathbb S^3$, the space $\Omega$ cannot be homemorphic to the real projective plane. So $\Omega$ is homeomorphic to $\frac{\mathbb S^1\times [0,1]}{\mathbb S^1\times\{0,1\}}$ and consequently, any neighborhood of $a$ in $\Omega$ cannot be topologically embedded in $\mathbb C$, because it contains the cone of two disjoint circles. Since $\Omega \subset A$, then $A$ is not locally embedded in $\mathbb C$.  \qed

\noindent \emph{Proof of Theorem~\ref{teoBom}.} As we have said, the proof is by induction on the dimension $n$. 

\smallskip
\noindent \emph{Case} $n=2$. Suppose $K\subset \mathbb{C}^2$ is a bombon, and let 
$\mathfrak{E}$ be the  ellipsoid of minimal volume containing $K$ (see Theorem~\ref{MiCE}). We may assume without loss of generality that $\mathfrak{E}$ is the unit ball of $\mathbb{C}^2$; and we must prove that $K = \mathfrak{E}$, which is equivalent to 
$$K\cap \partial \mathfrak{E}=K\cap \mathbb S^3=\mathbb S^3\,.$$

By Lemma \ref{lemclo},  $K\cap \mathbb S^3$ is a linearly closed subset of $LS^3$ and by Lemma \ref{prop:lin_cl}, we have to prove that $K\cap \mathbb S^3$ is neither a point nor an line. It clearly cannot be a point by the minimality condition. So, we are left to prove that this condition also rules out the line. We give an \emph{ad hoc} proof of this fact, without using John's theory on minimal ellipsoids (see for example \cite{J} or \cite{ball}).

Suppose there is a  line $L$, such that $K\cap \mathbb S^3=L\cap \mathbb S^3$.
First, we consider the case in which $L$ passes through the origin; without loss of generality, suppose
$L=\mathbb C\times\{0\}$. 

Consider the following family of ellipsoids for small $\epsilon\geq 0$, (see Fig. 1.a)
$$\mathfrak{E}_\epsilon= \{(x,y)\in \mathbb{C}^2\mid \frac{|x|^2}{ (1+\epsilon)^2}+ \frac{|y|^2}{ (1-\epsilon)^2}\leq 1\}\,.$$ 
For $\epsilon >0$, the volume of $\mathfrak{E}_\epsilon$ is less than the volume of the unit ball $\mathfrak{E}=\mathfrak{E}_0$; in fact, their volume ratio is $(1-\epsilon^2)^2$. Then, all the sets  $X_\epsilon=\mathfrak{E}\setminus \mathfrak{E}_\epsilon$ intersect $K$ because of the minimality condition of $\mathfrak{E}$. Furthermore, the intersection and limit of $X_\epsilon$, when $\epsilon \to 0$, is the solid torus $T=\{\,(x,y)\in\mathbb S^3\, |\, |y|\geq1/\sqrt{2}\,\}$. So that the compacity of $K$ implies that there exists a point in $K\cap T\subset K\cap \mathbb S^3$, which is far from $L$.

\begin{figure}[h!]
          \includegraphics[width=13cm]{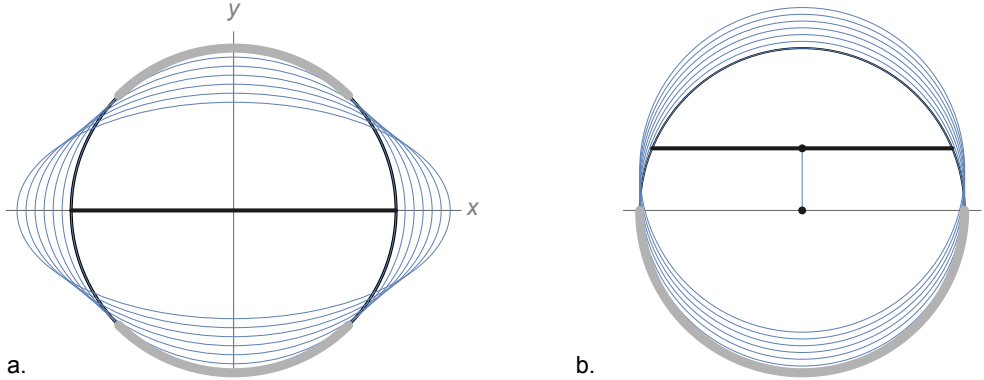}    
          \caption{\small {\bf a)} Section of the family of ellipsoids $\mathfrak{E}_\epsilon$, with a real plane having real lines in the coordinate axis.  {\bf b)} Section of the family of balls $\mathfrak{E}+\epsilon c$ with a real plane containing $c$ and the origin.}  
\end{figure}

If $L$ does not contain the origin, let $c$ be its closest point to the origin and consider the family of balls $\mathfrak{E}+\epsilon c$, $\epsilon\geq 0$ (see Fig.1.b). As before, the sets $X_\epsilon=\mathfrak{E}\setminus (\mathfrak{E}+\epsilon c)$ intersect $K$ for all $\epsilon>0$, because of the uniqueness  of the minimal volume ellipsoid containing $K$ (Theorem~\ref{MiCE}). Therefore, $K$ has a point in the limit of $X_\epsilon$, as $\epsilon \to 0$, which is the closed half-hemisphere of $\mathbb S^3$ opposite $c$. This contradicts that $K\cap \mathbb S^3=L\cap \mathbb S^3$, and completes the proof for $n=2$.

\noindent\emph{General case.} Suppose the theorem is true for $n-1$, we shall prove it for $n\geq3$.  By induction, for every hyperplane $H$, the section $H\cap K$ is either empty, a single point or an ellipsoid. This implies that every section of $K$ is symmetric and consequently, by Theorem~\ref{lemBB}, $K$ is symmetric.  

We may assume without loss of generality that the center of $K$ is the origin.  Now note that every hyperplane section of $K$ through the origin is an ellipsoid.
By Proposition \ref{prop:Celipsoo}, this implies that $K$ is an ellipsoid.  \qed

\section{Complex ellipsoids by means of sections and of projections}

The main purpose of this section is to prove that all  sections through a point $p_0$  of a convex body are  ellipsoids only when  the body is an ellipsoid.    The corresponding $3$-dimensional theorem for sections of real ellipsoids, as we already mentioned, is due to Kubota in \cite{Ku} and was independently discovered by Auerbach, Mazur and Ulam, \cite{AMU}. In its $n$-dimensional real version the result was originally proved by Busemann for $p_0\in  int(K)$ \cite{Bus}, and was extended to any $p_0\in \mathbb{E}^n$ by Burton, \cite{Bur1}. We also prove
the corresponding theorem for orthogonal projections. All orthogonal projections of a convex body onto affine $k$-planes are ellipsoids only when the body is an ellipsoid.  For real ellipsoids this result was observed by Blaschke and Hessenberg \cite{BH} without proof and proved by S\"uss \cite{Su0} in dimension $3$.  Later, Chakerian in \cite{Ch5}, gave a proof for the case of $n$ dimensions.

\begin{theorem}\label{teoksec}
Let $k$ and $n$ be such that $2\leq k < n$.  A convex body in $\mathbb C^n$ all whose  $k$-plane sections through a point $p_0$  are ellipsoids is an ellipsoid. 
\end{theorem}

\proof  It is easy to see that the general case $2\leq k < n$ of this theorem follows from the case $k=n-1$.  
Let $K\subset \mathbb{C}^n$ be a convex body, $n\geq 3$, and suppose  that all  hyperplane sections of $K$ through $p_0$ are ellipsoids.
Let us first prove that $K$ is a real ellipsoid.  Using Theorem 1.12.4 of \cite{MMO}, it is enough to prove that all real plane sections of $K$  are real ellipses.  Let $P$ be a real $2$-dimensional plane through $p_0$ and let $H$ be a  hyperplane through $p_0$ such that $P\subset H\subset \mathbb{C}^n$. This is possible because $n\geq 3$. By hypothesis, $H\cap K$ is an ellipsoid.  Then, $P\cap K$  is a real ellipse.  This implies that $K$ is a real ellipsoid. 

Since every  hyperplane section of $K$ through $p_0$ is an ellipsoid, we have that  
all  line sections of $K$ through $p_0$ are disks. Remember that any two parallel sections of a real ellipsoid are homothetic.  Hence, every complex line section, not only the ones through $p_0$, are disks and therefore $K$ is a bombon.  By Theorem \ref{teoBom}, $K$ is an ellipsoid. \qed

\medskip
An equivalent result is true for orthogonal projections.  

\begin{theorem}
Let $k$ and $n$ be such that $2\leq k < n$. A convex body in $\mathbb C^n$ all whose orthogonal projections onto $k$-dimensional planes  
are  ellipsoids is an  ellipsoid.
\end{theorem}

\proof Suppose the orthogonal projection of the convex set $K\subset \mathbb C^n$ onto every $k$-dimensional subspace is an ellipsoid, with $2\leq k < n$.
Since  ellipsoids are  symmetric, by Theorem \ref{teoSim}, $K$ is symmetric. By Lemma \ref{symellii}, it now suffices to prove that  $K$ is a real ellipsoid.

Let $P$ be a $2$-dimensional real plane and let $H$ be a complex $k$-subspace such that $P\subset H\subset \mathbb{C}^n$; this is possible because $k\geq 2$. By hypothesis, the orthogonal projection of $K$ onto $H$ is an ellipsoid.  Then, the orthogonal projection of $K$ onto $P$ is a real ellipse.  By Theorem 2.12.5 of \cite {MMO}, $K$ is a real ellipsoid.
 \qed

Our next purpose is to prove that  ellipsoids are the most symmetric bodies in complex space, in the sense that only ellipsoids have all their hyperplane sections symmetric.   For real ellipsoids, this result was first proved  by H. Brunn \cite{Bru} under the hypothesis of regularity and in general by G. R. Burton 
\cite{B2}.  It is interesting to note the use of  topology in the proof of the corresponding complex  result.

\medskip 

\begin{theorem}\label{teocs}
A convex body all whose  hyperplane sections are  symmetric is an ellipsoid. 
\end{theorem}

\proof  Let $K$ be a convex body all whose  hyperplane sections are  symmetric and let  $L\subset \mathbb{C}^{n-1}$ be a  line such that $L\cap K\not=\emptyset$. We shall prove that $L\cap K$ is  a disk and the theorem follows from Theorem~\ref{teoBom}. 
For that purpose, 
it will be enough to prove that there exists a hyperplane $H$ containing $L$,  such that the center of $H\cap K$ lies in $L$.  

Suppose, without loss of generality, that
$L=\{(0, \dots,0, z)\mid z\in \mathbb{C}\}$, and that the origin is in $L\cap K$.
For every  $(n-2)$-plane  through the origin  $\Gamma\subset \mathbb{C}^{n-1}$, we have that 
$\Gamma\cap K\not=\emptyset$. Let $\Gamma^\prime$ be the hyperplane of $\mathbb{C}^n$ generated by $\Gamma$ and $L$, and let $x_\Gamma$ be the center of symmetry of $\Gamma^\prime\cap K$. If $x_\Gamma \in L$,
there is nothing to prove; so we may assume, without loss of generality, that for every $(n-2)$-plane of $\mathbb{C}^{n-1}$  through the origin  $\Gamma$, we have that $x_\Gamma\notin L$.

Let $\pi: \mathbb{C}^{n}\to \mathbb{C}^{n-1}$ be the projection  onto the first $n-1$ coordinates. The choice $\Gamma\to \pi(x_\Gamma)\in\Gamma$ is a continuous assignment of a non-zero vector in $\Gamma$, for each hyperplane $\Gamma$ of $\mathbb{C}^{n-1}$. This is a contradiction to the well known fact that the canonical vector bundle of hyperplanes through the origin of $\mathbb{C}^{n-1}$ does not admit a non-zero section. See Steenrod's Book \cite{St}.

The continuity of $\Gamma\to \pi(x_\Gamma)$ follows from the following fact: suppose $H_i, i\geq 1$, is a sequence of hyperplanes in $\mathbb{C}^{n}$ with the property that 
$H_0=\lim 
(H_i)$ (where limits are taken as ${i\to\infty}$). If $H_i\cap K$ is symmetric with centre at $c_i$ for $i\geq 1$, then  $H_0\cap K$ is  symmetric with center at ${c_0=\lim c_i}$.
To see this, let $\Gamma_i=H_i-c_i$ and $K_i=(H_i\cap K)-c_i$, for $i\geq 0$.  Then, $\Gamma_0=\lim \Gamma_i$,  $K_0=\lim K_i$ (using the Hausdorff metric) and $K_i$ is symmetric with center at the origin for $i\geq 1$. This implies for every 
$\xi\in\mathbb{S}^1$ and $i\geq 1$,  that $\xi  K_i= K_i$. Therefore, 
$$\xi K_0=\xi(\lim K_i)=\lim \xi K_i =\lim  K_i= K_0$$ 
which implies that $H_0\cap K$ is symmetric with center at $c_0$.   \qed

\section{Convex bodies all whose sections are complex affinely equivalent}\label{Sect_Banach}

In 1932, Stephan Banach conjectured that if all $n$-dimensional subspaces of a Banach Space $V$ are isometric, $n>1$, then $V$ must be a Hilbert space. Of course, the Banach Space $V$ can be a linear space over the reals or the complex numbers.  If $V$ is a complex Banach space, this conjecture is equivalent to the following characterization of the ellipsoid. \emph{ A symmetric convex body $B\subset \mathbb C^{n+1}$ all whose  hyperplane sections through the origin are linearly equivalent  is an ellipsoid.}   This problem has a positive answer when $n$ is even (Gromov \cite{G}) and when $n\equiv 1$ mod $4$ (Bracho and Montejano \cite{BM}). The purpose of this section is to give a brief summary of the ideas and techniques used in the proof of the following theorem.

 \begin{theorem}\label{teocmain}
If all hyperplane sections through the origin of a symmetric convex body $B\subset \mathbb C^{n+1}$ are linearly equivalent, $n\equiv 0,1,2$ mod $4$, 
 then the convex body $B$ is an ellipsoid.
\end{theorem}

Given a  symmetric convex body  $K\subset \mathbb{C}^n$, let 
$$G_K:=\{\,g\in GL_n(\mathbb{C})\,|\, g(K)=K, \det g>0\,\}$$ 
be the  {\em group of linear isomorphisms of $K$ with positive real determinant}.   By Theorem \ref{MiCE}, there exists an ellipsoid of minimal volume containing $K$  centered at the origin. Suppose now that this minimal ellipsoid is the $(2n-1)$-dimensional unit ball, then every  $g\in G_K$ is actually an element of $SU_n$, because it fixes the unit ball, so in this case, $G_K=\{\,g\in SU_n\mid g(K)=K\,\}\subset SU_n$.

Let $B\subset \mathbb{C}^{n+1}$, $n\geq2$, be a symmetric convex body with center at the origin all whose hyperplane sections through the origin are  linearly equivalent.  By the above, it is not difficult to prove that  there exists a symmetric convex body $K\subset \mathbb{C}^n$ with center at the origin and with the property that every hyperplane section of $B$ is  linearly equivalent to $K$ and moreover $G_K\subset SU_n$. 
Indeed, using Topology of compact Lie groups, in particular, the notion of reduction of the structure group of a fiber bundle (see \cite{St} and \cite{L}), 
it is possible to prove that if $n$ is even, then $G_K=SU_{n}$ and if 
$n\equiv 1 $ mod $4$, then $SU_{n-1}\subset G_K$. For the details see \cite{BM}.  
 Note that if $n$ is even, due to the fact that $SU(n)$ is transitive, $G_K=SU_{n}$ implies that $K$ is a ball and hence that all hyperplane sections of $B$ are ellipsoids.  Consequently, by Theorem~\ref{teocs}, $B$ must be an ellipsoid. But  to understand the convex geometry of the consequences of the fact that $SU_{n-1}\subset G_K$, 
 when $n\equiv 1 $ mod $4$,  we need the following definition:
  
A {\em  body of revolution} is a symmetric convex body $K\subset \mathbb{C}^n$ for which there exists a \linebreak$1$-dimensional subspace $L$ of $\mathbb{C}^n$, called its {\em axis of revolution}, such that for every  affine hyperplane $H$ orthogonal to $L$, we have that $H\cap K$ is either empty, a single point or a $(2n-2)$-dimensional ball centered at $H\cap L$.  Of course, $K$ is a convex body of revolution if and only if $SU_{n-1}\subset G_K$.

With this in mind,  using topology of compact Lie groups, it is possible to prove  the following:

\begin{proposition}\label{thmcCkey} 
Let $B\subset \mathbb{C}^{n+1}$, $n\equiv 1$ mod 4,  $n\geq 5$, be a symmetric convex body with center at the origin all whose  hyperplane sections through the origin are  linearly equivalent. Then, there exists a  body of revolution $K\subset \mathbb{C}^n$ with center at the origin and with the property that every  hyperplane section of $B$ through the origin is linearly equivalent to $K$.
\end{proposition}

The next step is to use convex geometry to prove that if all hyperplane sections of a symmetric body are linearly equivalent to a  body of revolution then at least one of the sections must be an ellipsoid. As we will see, the proof of this result  is an interesting combination of ideas of algebraic topology and convex geometry.  We need first to understand the geometry of those  symmetric convex bodies which are linearly equivalent to a body of revolution.  From now on, let us call  those bodies \emph { linear  bodies of revolution}. Thus, a {\em linear body of revolution}
comes equipped with an {\em axis of revolution}, $L$, which is a line, and a {\em hyperplane of revolution}, $H$, which is a complementary hyperplane (but not necessarily orthogonal) to $L$, and it satisfies that all its sections with affine hyperplanes $H^\prime$  parallel to $H$ are either empty, a point or an ellipsoid centered at $L$ and homothetic to the ellipsoid $H\cap K$.

A non-spherical real body of revolution has only one real axis of revolution. The reason is that between $SO_{n-1}$ and $SO_n$ there is no connected subgroup. The analogy over the complex numbers is false, that is, between $SU_{n-1}$ and $SU_n$ there are several connected subgroups, however, as we will see next, a non-spherical  body of revolution has also a single axis of revolution.

 \begin{lemma}\label{lema:ctwo}
A linear  body of revolution $K\subset \mathbb{C}^n$, $n\geq 3$, admitting two different  hyperplanes of revolution, is an ellipsoid.
\end{lemma}

\proof By Theorem \ref{MiCE}, let $\mathfrak{E}$ be the unique ellipsoid of minimal volume centered at the origin containing $K$ and we may suppose, without loss of generality, that $\mathfrak{ E}$ is the unit ball. Since every symmetry of $K$ is a symmetry of the unit ball, our hypothesis now implies that $K$ is a  body of revolution with two different axes of revolution.  
Let $L_1$ and $L_2$ be two different  lines and let $G_1$ and $G_2$ be the complex rotation groups around the axis $L_1$ and $L_2$, respectively; they are both conjugate to $SU(n-1)$.
Suppose $G$ is a compact subgroup of $SU_n$ that contains both $G_1$ and $G_2$. We shall prove that the action of $G$ in $\mathbb{S}^{2n-1}$ is transitive. If this is so, and both $L_1$ and $L_2$ are axes of revolution of $K$, then $G^0_K=\{g\in SU_n\mid g(K)=K\}$, which is compact because $K$ is a compact convex body, would act transitively on $\partial K$ and $K$ would be a ball.

 Let $P$ be the plane generated by $L_1$ and $L_2$ and let $\pi_1, \pi_2$ and $\pi_0$ be the orthogonal projections onto $L_1, L_2$ and $P$, respectively. Furthermore, let $D=P\cap  int(\mathfrak{B})$, where $\mathfrak{B}$ is the unit ball of $\mathbb{C}^n$.  
Consider the set
$$U=\pi^{-1}_0(D)\cap\mathbb{S}^{2n-1}\,.$$
Note that $U$ is an open connected dense subset of $\mathbb{S}^{2n-1}$ because $\mathbb{S}^{2n-1}\setminus U=P\cap \mathbb{S}^{2n-1}$ is a $3$-sphere contained in $\mathbb{S}^{2n-1}$, and since $n\geq 3$, its (topological) codimension is at least 2. 

Let $x\in U$. Our purpose is to construct an open neighborhood $W$ of $x$ in $U$ such that $W$ is contained in the orbit $G \cdot x$ of $x$ under the action of $G$ in $\mathbb{S}^{2n-1}$. This will be enough to prove the lemma because $U$ is a connected open dense subset of $\mathbb{S}^{2n-1}$. 

Let $H_1=\pi_1^{-1}(\pi_1(x))$. It is the  affine hyperplane orthogonal to $L_1$ and passing through $x$, so that $G_1 \cdot x=H_1\cap\mathbb{S}^{2n-1}$. Let $W_1=H_1\cap D$. It is an open disk in an affine line parallel to the line $L_1^\perp\cap P$, and observe that restricted to this affine line ($H_1\cap P$), the map $\pi_2$ is an affine isomorphism onto $L_2$ because $L_1 \not= L_2$. So that $W_2 = \pi_2(W_1)$ is an open subset of $L_2 \cap D$ that contains $\pi_2(x)$.

Let $W=\pi_2^{-1}(W_2)\cap U$. It is an open neighborhood of $x$ in $U$. We are left to prove that $W$ is contained in the orbit $G \cdot  x$. 

Given $y\in W$, let $H_2=\pi_2^{-1}(\pi_2(y))$, so that $G_2 \cdot  y=H_2\cap\mathbb{S}^{2n-1}$. Consider the affine subspace $\Gamma=H_1 \cap H_2$ of dimension $n-2 >0$. By construction, $H_2$ intersects $W_1$ in a point, so that $\Gamma$ touches the interior of the unit ball $\mathfrak{B}$. Therefore, $\Gamma \cap \mathbb{S}^{2n-1} = (G_1 \cdot x)\cap(G_2 \cdot  y)$ is not empty. This implies that $G \cdot x=G \cdot  y$, so that $y\in G \cdot x$, and hence $W\subset G \cdot x$. 
\qed

\begin{lemma}\label{lemCsubset} 
Every  hyperplane section $\Gamma \cap K$ of a linear body of revolution ${K\subset \mathbb{C}^{n}}$, $n\geq 3$,  is a linear body of revolution. Furthermore, if $H$ is the  hyperplane of revolution of $K$, then either $\Gamma=H$ or $\Gamma\cap H$ is a hyperplane of revolution of  $\Gamma \cap K$.  
\end{lemma}

\proof Without loss of generality, we may assume  that $K$ is a  body of revolution; that is, if its axis of revolution is the  line $L$, then $H=L^\perp$ is the corresponding  hyperplane of revolution and we have that $H\cap K$ is a ball centered at the origin.  

Assume $\Gamma\not=H$. We will prove that  $K_1=\Gamma\cap K$ is a body of revolution in $\Gamma$ with  hyperplane of revolution $H_1=\Gamma\cap H$. Let $L_1$ be the  line orthogonal to $H_1$ in $\Gamma$; it will be the axis of $K_1$. 

Given $H_1^\prime\subset\Gamma$ parallel to $H_1$, we have to consider the intersection of $H_1^\prime\cap K_1=H_1^\prime\cap K$.
 Let $H^\prime$ be the affine  hyperplane of $\mathbb{C}^{n}$ parallel to $H$ that contains $H_1^\prime$. By hypothesis, we have that $H^\prime\cap K$ is either empty, a point or a ball (in $H^\prime$) centered at $L$. Therefore, its intersection with $H_1^\prime$ (a hyperplane of $H^\prime$) is either empty, a point or a ball. By construction, in the two last cases the point or the center of the ball lies in $L_1$; indeed, the plane generated by $L$ and $L_1$ is orthogonal to $H_1$. Therefore, $K_1\subset\Gamma$ is a  body of revolution as we wished.
\qed

\begin{lemma} \label{lemCE}
Let $K\subset \mathbb{C}^n$ be a linear body of revolution with axis of revolution $L$, $n\geq3$.  Suppose $\Gamma\subset \mathbb{C}^n$ is a hyperplane containing $L$ for which $\Gamma\cap K$ is an ellipsoid. 
Then $K$ is an ellipsoid. 
\end{lemma}

\proof First, we may assume that $K$ is a  body of revolution with axis of revolution $L$, hyperplane of revolution $H=L^\perp$ and such that $H\cap K$ is the unit ball in $H$.  By hypothesis and Lemma \ref{lemCsubset}, $\Gamma\cap K$ is an ellipsoid and a  body of revolution with axis of revolution $L$. 
Using a linear map which is the identity on $H$ and a dilatation on $L$, we may assume $\Gamma\cap K$ is a unit ball centered at the origin; so that $\Gamma\cap K=\Gamma\cap\mathfrak{B}$, 
where $\mathfrak{B}\subset \mathbb{C}^n$ is the unit ball. Our purpose is to prove that $K=\mathfrak{B}$ to conclude the proof.  

For every affine hyperplane $H^\prime$ parallel to $H$ that touches the interior of $K$, we have that both $H^\prime\cap K$ and  $H^\prime\cap \mathfrak{B}$ are concentric balls. Furthermore, they have the same radius because their boundaries have non empty intersection (in $\Gamma$). Consequently, $H^\prime\cap K=H^\prime\cap \mathfrak{B}$ and hence $K=\mathfrak{B}$, as we wished.  
\qed

From now on, let   $B\subset \mathbb{C}^{n+1}$  be a symmetric convex  body, $n$ odd, $n\geq 5$, all of whose  hyperplane sections are  non elliptical, linear bodies of revolution. First of all note that the assignment of the axis of revolution in every hyperplane is continuous. 

For every line $\ell\subset \mathbb{C}^{n+1}$  denote by $\ell^\perp$ the  hyperplane subspace of $\mathbb{C}^{n+1}$ orthogonal to $\ell$.  Furthermore, by Lemma \ref{lema:ctwo}, denote by $L_\ell$ the unique  axis of revolution of $\ell^\perp \cap B$ 
and by $H_\ell$ the corresponding  $(n-1)$-dimensional subspace of revolution of $\ell^\perp\cap B$. 
Note that the line $L_\ell$ contains the origin because the center of the ellipsoid $H_\ell \cap B$ is
 the origin.

\begin{lemma}\label{lemprin} Let $B\subset \mathbb{C}^{n+1}$, $n$ odd, $n\geq 5$,  be a symmetric convex  body, all whose hyperplane sections are  non elliptical, linear bodies of revolution.
Suppose  $\ell_1$ and $\ell_2$ are two different  lines with the property that $L_{\ell_2}\subset \ell_1^\perp$. Then
 $$ H_{\ell_2}\cap\ell_1^\perp=H_{\ell_1}\cap\ell_2^\perp=H_{\ell_2}\cap H_{\ell_1}.$$
 \end{lemma}
 
 \proof Consider $\ell_1^\perp\cap\ell_2^\perp$, the $(n-1)$-dimensional subspace of $\ell_2^\perp$. By hypothesis, 
 $\ell_2^\perp \cap B$ is a non  elliptical, linear body of revolution
  with axis of revolution $L_{\ell_2}$. Therefore, since $L_{\ell_2}\subset \ell_1^\perp\cap\ell_2^\perp$,
  we have that $\ell_1^\perp\cap\ell_2^\perp\cap B$
  is a  linear body of revolution with corresponding  axis of revolution $L_{\ell_2}.$
  Furthermore, by Lemma \ref{lemCE}, $\ell_2^\perp \cap B$ is not an ellipsoid. Moreover, the unique  hyperplane of revolution of 
  $\ell_1^\perp\cap\ell_2^\perp \cap B$ is $\ell_1^\perp \cap H_{\ell_2}.$
 
 On the other hand, $\ell_1^\perp\cap\ell_2^\perp$ is an 
 $(n-1)$-dimensional subspace of $\ell_1^\perp.$ 
 Note that $\ell_1^\perp\cap\ell_2^\perp\not=H_{\ell_1}$, otherwise $\ell_1^\perp\cap\ell_2^\perp\cap B=H_{\ell_1}\cap B$
  would be an ellipsoid contradicting our previous assumption. Since $\ell_1^\perp\cap B$ is a non elliptical, linear body of revolution and $\ell_1^\perp\cap\ell_2^\perp\not= H_{\ell_1}$, then by Lemma \ref{lemCsubset}, 
  $\ell_1^\perp\cap\ell_2^\perp\cap B$ is a
 non elliptical, linear body of revolution with corresponding  subspace of revolution $H_{\ell_1}\cap \ell_2^\perp.$ Consequently, by Lemma \ref{lema:ctwo}, since $n-1\geq3$, we have that $ H_{\ell_2}\cap\ell_1^\perp=H_{\ell_1}\cap\ell_2^\perp$.
 \qed
 
 \smallskip
 We are ready to prove that 
 \smallskip
 
\begin{theorem}\label{thm:Celip} Let $B\subset \mathbb{C}^{n+1}$  be  a symmetric convex  body, $n$ odd, $n\geq 5$. Suppose all hyperplane sections of $B$ are linear bodies of revolution. Then, one of the hyperplane sections is an ellipsoid. 
\end{theorem}

\proof The assignment $\ell \to L_\ell\subset \ell^\perp$ is a continuous function on $\ell$. 
 We start proving that there are two different  lines $\ell_1$ and $\ell_2$ such that
$$ L_{\ell_2}\subset H_{\ell_1}.$$
For that purpose, let $L$ be a line in $H_{\ell_1}$ and denote by $L^\perp$ the  hyperplane of $\mathbb{C}^{n+1}$ orthogonal to $L$.
Denote by  $\Pi:\mathbb{C}^{n+1} \to L^\perp$ the orthogonal projection along $L$. Of course, $\Pi$ is a linear map over the complex numbers  $\mathbb{C}$. 

For every  line $\ell\subset L^\perp$, consider $L_\ell\subset \ell^\perp$. Suppose $L_\ell\not=L$, for every $\ell\subset L^\perp$. If this is so, $\Pi(L_\ell)$ is a line through the origin in $L^\perp$, ortogonal to $\ell$. Furthermore, the assignment $\ell\to \Pi(L_\ell)$ is a continuous function on $\ell$.  This is a contradiction to the well known fact that the canonical vector bundle of hyperplanes through the origin of $\mathbb{C}^{n-1}$ does not admit a 
 line subbundle, when $(n-1)$ is even.  See Steenrod's Book \cite{St}.
 This implies the existence of $\ell_2$ such that $L=L_{\ell_2}\subset H_{\ell_1}$.

Note now that this is a contradiction to Lemma \ref{lemprin} because clearly $L_{\ell_2}\subset \ell_1^\perp$, hence
$$L_{\ell_2}\subset H_{\ell_1}\cap\ell_2^\perp = H_{\ell_2}\cap\ell_1^\perp \subset H_{\ell_2}$$
which is impossible. 
\qed

Once we know at least one of the sections of our body $B$ has a hyperplane section through the origin which is an  ellipsoid, all hyperplane sections through the origin are  ellipsoids and hence, by Theorem \ref{teocs}, $B$ is an ellipsoid. This concludes the proof of Theorem \ref{teocmain}.

The following, and final, result is again a characterization of ellipsoids by means of their  line sections and following the spirit of the above results.  As the proof is shorter and simpler than the proof of Proposition \ref{thmcCkey}, we include it in its entirety  so that the reader has a more complete idea of the topological ideas developed in this section.

\begin{theorem}\label{entenado}
 A convex body $K$ is an ellipsoid  if and only if any two  non-trivial line sections of $K$ are 
affinely equivalent. 
\end{theorem}

\proof  In view of  Theorem~\ref{teoBom}, it will suffice to prove that every non-trivial  line section of $K$ is a disk. 

By hypothesis, there is a fixed convex set $C\subset \mathbb{C}$, such that every non-trivial  line section of $K$  is  affinely equivalent to $C$. Without loss of generality, we may assume that $\mathbb{S}^1$ is the circumcircle of $C$, so that
$G_C$, \emph{the compact group of   affine equivalences of $\mathbb{C}$ that fix $C$}, is a subgroup of $U(1)$. To see that every non-trivial line section of $K$ is a disk, we will prove that $C$ is the unit disk, which is clearly equivalent to  proving that $G_C=U(1)$.

Consider a point in the interior of $K$ and a 2-dimensional  affine subspace through it. To fix ideas, we may assume that the origin is in the interior of $K$, and think of $\mathbb{C}^2\subset\mathbb{C}^n$.  The space of all  lines through the origin in $\mathbb{C}^2$ is the Riemann sphere $\mathbb{S}^2$. Let 
$$\xi^2:\mathbb{\it E}^2\to \mathbb{S}^2$$
be the canonical vector bundle of $1$-dimensional  subspaces of $\mathbb{C}^2$. 
In view of the continuity of the  line sections of $K$, our hypothesis imply that 
the structure group of $\xi^2$ can be taken to be $G_C$. That is, if  $\xi^2$ is described by a characteristic map $\tau:\mathbb{S}^1\to U(1)$, it factors through $G_C$; (or if it is described by  transition functions and an open cover, the transition functions can be chosen to go to $G_C$).  See, e.g., \cite{St}. 

Since the only proper compact subgroups of $U(1)$ (which can be identified with $\mathbb{S}^1$) are finite, but the canonical  line bundle is non trivial, then $G_C=U(1)$. This completes the proof.
 \qed
 
\bigskip

\noindent {\bf Acknowledgments.} Luis Montejano acknowledges  support  from CONACyT under project 166306 and  from PAPIIT-UNAM under project IN112614.

\end{document}